\documentclass[12pt]{article}
\usepackage{amsfonts}
\usepackage{amsfonts,amssymb,amsmath,mathrsfs}
\usepackage{color,latexsym,amsfonts}
 \newcommand{\qed}{\hfill\rule{2mm}{3mm}\vspace{4mm}}

 \newtheorem{theorem}{Theorem}[section]
 \newtheorem{lemma}[theorem]{Lemma}
 \newtheorem{corollary}[theorem]{Corollary}
 \newtheorem{proposition}[theorem]{Proposition}
 \newtheorem{example}[theorem]{Example}
 \newtheorem{remark}{\noindent\mbox{Remark}}[section]

 \def\beqlb{\begin{eqnarray}}\def\eeqlb{\end{eqnarray}}
 \def\beqnn{\begin{eqnarray*}}\def\eeqnn{\end{eqnarray*}}
\def\red{  }

 \def\proof{\noindent{\it Proof.~~}}
 \def\qed{\hfill$\Box$\medskip}
 \def\no{\nonumber}
 \def\red{  }

\title{Light-tailed behavior of stationary distribution
for state-dependent random {  walks} on a
strip\thanks{\red{The project was partially supported by the
National Natural Science Foundation of China (Grant No. 11131003)
and by the Natural Sciences and Engineering Research Council of
Canada (Grant No. 315660).}}}

\author{{  Wenming Hong\footnote{School of Mathematical Sciences
\& Laboratory of Mathematics and Complex Systems, Beijing Normal
University, Beijing 100875, P.R. China. Email: wmhong@bnu.edu.cn},
\ \ \
 Meijuan Zhang\footnote{ School of Mathematical Sciences
\& Laboratory of Mathematics and Complex Systems, Beijing Normal
University, Beijing 100875, P.R. China. Email:
zhangmeijuan@mail.bnu.edu.cn} \ \ \ and  \ \ \ Yiqiang Q.
Zhao}\footnote{School of Mathematics and Statistics, Carleton
University, Ottawa, Ontario, Canada K1S 5B6. Email:
zhao@math.carleton.ca}}
\date{}

\begin{document}
\maketitle

\begin{center}
\begin{minipage}[c]{15cm}

\mbox{}\textbf{Abstract:}\quad  {  In this paper, we
consider} the state-dependent reflecting random walk on a
half-strip.{We provide} explicit criteria for
(positive) recurrence, and {  an explicit} expression for
the stationary distribution. {As a consequence,} the
light-tailed behavior of the stationary distribution is proved
{  under appropriate conditions. The key idea of the
method employed here is the decomposition} of the trajectory of
the random walk and the main {  tool is} the intrinsic
branching structure {  buried in} the random walk on a
strip, {  which is different from the matrix-analytic
method.} \vspace{0.2cm}

\mbox{}\textbf{Keywords:}\quad random walk on a strip, stationary distribution, light-tailed behavior, branching process, recurrence, state-dependent.
\vspace{0.2cm}

\mbox{}\textbf{Mathematics Subject Classification}: Primary 60K37;
secondary 60J85.
\end{minipage}
\end{center}

\section{Introduction and Main results}

Let $d\geq1$ be any integer and denote
$\mathscr{D}=\{1,2,\cdots,d\}$. {  We consider the}
reflecting space-inhomogeneous and state-dependent reflecting
random walk on a half-strip {   $S=\{0, 1, 2, \ldots \}
\times\mathscr{D}$.} This model is often referred to as the
state-dependent quasi-birth-and-death (QBD) process in queueing
theory. {  Studies on the state-dependent QBD process has
been centered at its stationary distribution such as properties of
the rate matrices, efficient algorithms for computations, often
through the matrix-analytic approach or the censoring techniques
(e.g. \cite{[LR99]}, \cite{[BT]} and \cite{[ZLB03]}). In this
paper, we propose a different method to decompose the trajectory
of the random walk on the strip using the intrinsic branching
structure ( \cite{[HM12]}), through which, we provide criteria for
(positive) recurrence, obtain an expression for the stationary
distribution of the walk, and characterize the exponential tail
asymptotic behavior of the stationary distribution for the walk.}

We consider the random walk $X_{n}=(\xi_{n},Y_{n})$,
{  $n =0, 1, \ldots$, on a half-strip, which is a Markov
chain with $\xi_{n}\in \{0, 1, 2, \ldots \}$, referred to as the
layer (or level),} and $Y_{n}\in\mathscr{D}$. Let
$L_{i}=\{(i,r);~r=1,2,\cdots,d\}$. {  Then, the
half-strip $S$} can be expressed as $S=\cup_{i=0}^{+\infty}L_{i}$.
The transition probability of {  the walk is given by}
\begin{equation}\label{tp}
\breve{ P}=\left(
         \begin{array}{cccccc}
               \mathbf{0}    & P_{0} & ~    & ~   & ~      &~ \\
               Q_{1}& R_{1} & P_{1}& ~   & ~      & ~\\
               ~    & Q_{2}& R_{2}& P_{2}& ~      & ~\\
               ~    &  ~   & Q_{3}& R_{3} & P_{3} & ~\\
               ~    &  ~   &  ~   & \ddots & \ddots &\ddots \\
             \end{array}
           \right),
\end{equation}
where $P_{0}$ {  is a $d\times d$ stochastic matric and
$\{(P_{n},Q_{n},R_{n})$, $n\in \mathbb{Z}^{+}=\{1, 2, \ldots\}\}$,
satisfies $(P_{n}+Q_{n}+R_{n})\mathbf{1}=\mathbf{1}$ with
$\mathbf{1}$ being a column vector of ones. For a matrix $A$, its
$(i,j)$th component is denoted by $A(i,j)$, $1\leq i, j\leq d$.
Following  \cite{[BG00]}, we assume the following conditions,
under which the process is irreducible:}

\begin{description}
\item[$C_{1}$.] $\log(1-\| R_{n}+P_{n}\|)^{-1}<\infty$
{  and} $\log(1-\| R_{n}+Q_{n}\|)^{-1}<\infty$;

\item[$C_{2}$.] {  For any $n$ and any $j$,
$\sum_{i=1}^{d}P_{n}(i,j)>0$ and} $\sum_{i=1}^{d}Q_{n}(i,j)>0$;

\item[$C_{3}$.] The layer $0$ is in one communication class.
\end{description}

{   To state the main results, we introduce the following
notation or definitions: $\mathbf{e}_{i}=(\underbrace{0,\ldots,0,
1,0, \ldots,0}_{i\text{th component}})$; $I_{A}$ is the indicator
function of the set $A$; $B'$ and $\mathbf{b}'$ are the transposes
of matrix $B$ and vector $\mathbf{b}$, respectively;
$\mathbf{1}=(1,1,\cdots,1)^{'}$; $P_{\mu_n}(\cdot)$ represents the
probability given that the random walk starts from layer $n$ with
the distribution $\mu_n$;  $P_{(n,i)}(\cdot)$ means the
probability given that the walk starts from the site (state)
$(n,i)$; $E_{\mu_n}(\cdot)$ and $E_{(n,i)}(\cdot)$ are similarly
defined.

Define the hitting times $T_{n}$ and $T_{n}^{+}$ by
$T_{n}=\inf\{t\geq 0:~X_{t}\in L_{n}\}$ and $T_{n}^{+}=\inf\{t\geq
1:~X_{t}\in L_{n}\}$, respectively.} Define $f_{x,y}^{(n)}$ and
$f_{x,y}$ respectively as
$f_{x,y}^{(n)}=P_{x}(T_{y}^{+}=n)=P_{x}(X_{n}\in L_{y},~X_{m}
\notin L_{y},~1\leq m<n)$ and
$f_{x,y}=\sum_{n=1}^{+\infty}f_{x,y}^{(n)}$. {  Also,
define}
\begin{equation} E_{x}(T_{y}^{+})=
\begin{cases}
\sum_{n=1}^{+\infty}nf_{x,y}^{(n)}& \mbox{if~} f_{x,y}=1,\\
+\infty & \mbox{if~} f_{x,y}<1.
\end{cases}\no
\end{equation}

\noindent\textbf{Definition 1.} The layer $y$ {   is
(layer) recurrent} if $f_{y,y}=1$, otherwise it is {
(layer) transient}. If $E_{y}(T_{y}^{+})<+\infty$, the layer $y$
is called {   (layer) positive recurrent.

Define recursively for
 $n \in \{0, 1, 2, \ldots \}$,
\begin{equation}
\zeta_{0}^{+}=P_{0} \quad \mbox{and} \quad
\zeta_{n}^{+}=(I-Q_{n}\zeta_{n-1}^{+}-R_{n})^{-1}P_{n},  \quad n
=1, 2, \ldots.
\end{equation}}
The existence of $(I-Q_{n}\zeta_{n-1}^{+}-R_{n})^{-1}$ is a
consequence of {  assumption $(C_2)$. Note that $P_{0}$
is stochastic, so is $\zeta_{n}^{+}$ (see \cite{[BG00]}).
 {  Define also $\mathbf{u}_{0}^{+}=\mathbf{1}$, and}
for $n\geq1$,
\begin{equation}\label{a1}
A_{n}^{+}=(I-Q_{n}\zeta_{n-1}^{+}-R_{n})^{-1}Q_{n}\quad\mbox{and}\quad
\mathbf{u}_{n}^{+}=(I-Q_{n}\zeta_{n-1}^{+}-R_{n})^{-1}\mathbf{1}.
\end{equation}}

\begin{remark}{  $\zeta_{n}^{+}(i,j)$ has an interpretation as the
probability of the random walk starting from $(n,i)$ and with the
reflection at layer $0$ reaches layer $n+1$ at point
$(n+1,j)$, often referred to as the exit probability and denoted as $\eta_{n}$ as well in our model. Also, $A_{n}^{+}(i,j)$ can be interpreted as
the expected number  } {     of steps from   $(n,i)$ to $(n-1,j)$ caused by a step from layer $n+1$ to $(n,i)$ (see ($\ref{ma}$), i.e., the mean offspring of the ``father" step from layer $n+1$ to $(n,i)$).}
{  It is worthwhile to
mention that the rate matrix and the fundamental period matrix are
key probabilistic quantities in studying the level-independent QBD
process, They are generalized into two matrix sequences $R^+_n$
and $G^-_n$, respectively, when the method is used to study the
level-dependent QBD process (e.g., \cite{[LR99]} and \cite{[ZLB03]}). Their dual versions $R^-_n$
and $G^+_n$ also play an important role in the study using the matrix-analytic method (e.g. \cite{Miyazawa-Zhao:04}).
The matrix $\zeta_{n}^{+}$ is the same as $G^+_n$, while $A_{n}^{+}$ is a unique quantity from the branching process method, which is
not a usual measure used in the matrix-analytic method.
 }
\end{remark}

{   The first group of results are conditions for recurrence and positive recurrence of the walk.}
\begin{theorem}\label{thm1}
{  For} the random walk starts from layer $0$ with an initial distribution $\mathbf{\mu}$,
define
\begin{equation}
\beta^{+}= \sum_{k=0}^{+\infty}\mathbf{\mu}_{k}A_{k}^{+}A_{k-1}^{+}\cdots A_{1}^{+}\mathbf{1},\no
\end{equation}
where $\mathbf{\mu}_{k}=
\mathbf{\mu}\zeta_{0}^{+}\zeta_{1}^{+}\cdots\zeta_{k-1}^{+}$. Then the random walk is recurrent if and only if $\beta^{+}=\infty$.
\end{theorem}

\begin{remark}
Actually, $\beta^{+}$ in the Theorem is the expectation number of the visiting times by the random walk at   layer $0$ , which can be calculated by the means of the  intrinsic branching structure within the walk.
\end{remark}

To state the criteria for the positive recurrence, we need the {  ``exit probability"} from the other direction. {   Let $\mathbb{Z}=\{0, \pm 1, \pm 2, \ldots\}$ and $a \in \mathbb{Z}$.  For $n \leq a$, define recursively}
\begin{equation}
{  \zeta_{a,a}^{-}} = \rho \quad \mbox{and} \quad
\zeta_{n,a }^{-}=(I-P_{n}\zeta_{n+1, \  a}^{-}-R_{n})^{-1}Q_{n},~~ {  n<a,}
\end{equation}
where $\rho$ is stochastic, i.e.,  $\rho\mathbf{1}=\mathbf{1}.$ Then under condition $C$, the limit $\zeta_{n}^{-}=\lim_{a\to\infty}\zeta_{n,a}^{-}$ exist and satisfy the following equation (Theorem 1, \cite{[BG00]}),
\begin{equation}
\zeta_{n}^{-}=(I-P_{n}\zeta_{n+1}^{-}-R_{n})^{-1}Q_{n},~~n\in \mathbb{Z}.
\end{equation}
Define for $n\geq1$,
\begin{equation}\label{a2}
A_{n}^{-}=(I-P_{n}\zeta_{n+1}^{-}-R_{n})^{-1}P_{n}\quad\mbox{and}\quad \mathbf{u}_{n}^{-}=(I-P_{n}\zeta_{n+1}^{-}-R_{n})^{-1}\mathbf{1}.
\end{equation}

\begin{remark} {  $\zeta_{n}^{-}(i,j)$ has an interpretation as the
probability of the random walk starting from $(n,i)$ reaches layer $n-1$ at point
$(n-1,j)$, which is the same as $G^-_n$. Also, $A_{n}^{-}(i,j)$ is the same as $A_{n}^{+}(i,j)$ but from the other direction, can be interpreted as
the expected number } {     of steps from   $(n,i)$ to $(n+1,j)$ caused by a step from layer $n-1$ to $(n,i)$,}
{   which is unique to the branching process method.}
\end{remark}

{  For convenience, let $\mathbf{u}_{0}^{-}=1$.}
\begin{theorem}\label{thm2}
{  For the random walk starting from layer $0$ with an initial distribution $\mathbf{\mu}$, define}
\begin{equation}
\varrho^{+}_1= \mathbf{1}'P_{0}
\left(\sum_{k\geq1}A_{1}^{-}A_{2}^{-}\cdots A_{k-1}^{-}\mathbf{u}_{k}^{-}\right)+d. \no
\end{equation}
{  Then} the random walk is positive recurrent if and only if $\varrho^{+}_1<\infty.$
\end{theorem}

\begin{remark}\label{r1}
Actually, $\varrho^{+}= \mathbf{\mu}P_{0}
\left(\sum_{k\geq1}A_{1}^{-}A_{2}^{-}\cdots A_{k-1}^{-}\mathbf{u}_{k}^{-}\right)+\mathbf{\mu}\mathbf{1}$  is the expectation  of the first return  time of the random walk start at layer $0$ , which can be calculated {  by means of} the  intrinsic branching structure within the walk. Note that $\varrho^{+}\leq \varrho^{+}_1$, the criteria is independent of the initial distribution of the walk start at layer $0$.
\end{remark}

{  When the walk is state-independent, i.e., $(P_n, Q_n, R_n)=(P, Q, R)$ for $n>0$,  we} denote the walk as $\{\overline{X}_{n},~n\geq 0\}$, and have correspondingly
$$\zeta^{-}=(I-P\zeta^{-}-R)^{-1}Q,$$  and
\begin{equation}\label{a}
A^{-}=(I-P\zeta^{-}-R)^{-1}P ,~  \mathbf{u}^{-}=(I-P\zeta^{-}-R)^{-1}\mathbf{1}.
\end{equation}

\begin{corollary}\label{cor5}
Suppose that the random walk $\{\overline{X}_{n},~n\geq 0\}$ starts from layer $0$ with an initial distribution $\mathbf{\bar{\mu}}$. Then

\noindent (1) The random walk is positive recurrent if and only if
\begin{equation}
\bar{\varrho}^{+}_1= \mathbf{1}'P~(\sum_{k\geq1}~(A^{-})^{k-1}~\mathbf{u}^{-})
+d<\infty,\no
\end{equation}
\noindent (2) Denote the maximum eigenvalues of $A^{-}$ as $\lambda_{A^{-}}$. Then $\lambda_{A^{-}}<1$ whenever $\bar{\varrho}^{+}_1<\infty$.
\end{corollary}

{  We now state the main result for the stationary distribution. We assume that} the walk is positive recurrent and start from layer $0$ with {  a ``proper" distribution.}
 The so called ``proper" distribution is the ``censored measure", {  a terminology borrowed from queueing theory (e.g. \cite{[ZLB03]}). Define $\breve{P}_{1}$ by}
 \begin{equation}
\breve{P}_{1}=\left(
         \begin{array}{ccccc}
                 R_{1} & P_{1}& ~   & ~      & ~\\
                 Q_{2}& R_{2}& P_{2}& ~      & ~\\
                 ~   & Q_{3}& R_{3} & P_{3} & ~\\
                 ~   &  ~   & \ddots & \ddots &\ddots \\
             \end{array}
           \right).\no
\end{equation}
Let $S_{0}=L_{0}$ and {  let} $S_{1}=S / S_{0}$ be a partition of the state space $S$. {  Then} $\breve{P}$ can be partitioned according to $S_{0}$ and $S_{1}$ as
\begin{equation}
\breve{P}= \left(
         \begin{array}{cc}
                  \breve{P}^{0} & U \\
                  D             & \breve{P}_{1}\\
         \end{array}
           \right),\no
\end{equation}
where $\breve{P}^{0}={R}_{0}$,~$U=(P_{0},O,O,\ldots)$ and $D=(Q_{1},O,O,\ldots)^{T}$.

The censored matrix $\breve{P}^{(S_{0})}$ of $\breve{P}$ with the censoring set $S_{0}$ is defined by
\begin{equation}
\breve{P}^{(S_{0})}=\breve{P}^{0}+U\widehat{\breve{P}^{1}}D,\no
\end{equation}
where $\widehat{\breve{P}^{1}}=\sum_{k=0}^{+\infty}(\breve{P}^{1})^{k}$ is called the fundamental matrix of $\breve{P}^{1}$. $\breve{P}^{(S_{0})}$ is a $d\times d$ matric, and the censored matrix $\breve{P}^{(S_{0})}$ has a probabilistic interpretation: it is the probability {  that} the next state visited in $S_{0}$ is $j$, given that the process starts in state $i\in S_{0}$.

A measure $\mathbf{\breve{\mu}}_{0}$ which satisfies
\begin{equation}
\mathbf{\breve{\mu}}_{0}\breve{P}^{(S_{0})}=\mathbf{\breve{\mu}}_{0}
\end{equation}
is called as censored measure with censoring set $S_{0}$.

\begin{theorem}\label{thm3}
If $\varrho^{+}_1<\infty$ and the walk starts from layer $0$ with {  the} censored measure $\mathbf{\breve{\mu}}_{0}$, then the stationary distribution
$\{\mathbf{\nu}_{n}, ~ {  n =0, 1, 2, \ldots} \}$ exists and unique, which can be expressed explicitly as
\begin{equation}\label{vn}
\mathbf{\nu}_{n}
=\frac{\mathbf{\breve{\mu}}_{0}P_{0}A_{1}^{-}A_{2}^{-}\cdots A_{n-1}^{-}\widetilde{u}_{n}^{-}}{\mathbf{\breve{\mu}}_{0}P_{0}(\sum_{k \geq 1}A_{1}^{-}A_{2}^{-} \cdots A_{k-1}^{-}\mathbf{u}_{k}^{-})+\mathbf{\breve{\mu}}_{0}\mathbf{1}},\quad n>0;
\end{equation}
and \begin{equation}
\nu_{0}=\frac{\mathbf{\breve{\mu}}_{0}}{\mathbf{\breve{\mu}}_{0} P_{0}(\sum_{k \geq 1}A_{1}^{-}A_{2}^{-} \cdots A_{k-1}^{-}\mathbf{u}_{k}^{-})+\mathbf{\breve{\mu}}_{0}\mathbf{1}},
\end{equation}
where $\mathbf{\breve{\mu}}_{0} P_{0}(\sum_{k \geq 1}A_{1}^{-}A_{2}^{-} \cdots A_{k-1}^{-}\mathbf{u}_{k}^{-})+\mathbf{\breve{\mu}}_{0}\mathbf{1}<\varrho^{+}_1<\infty$ and $\widetilde{u}_{n}^{-}=(I-P_{n}\zeta_{n+1}^{-}-R_{n})^{-1}$.
\end{theorem}

{
\begin{remark}
We can show that the expression in Theorem~\ref{thm3} is consistent with the matrix-product form solution given by the matrix-analytic method:
\[
    \nu_n = \nu_0 R^+_1 R^+_2 \cdots R^+_n.
\]
To see it, we notice that $\tilde{u}^-_n = (I-P_n
\zeta^-_{n+1}-R_n)^{-1}$ is the entry $\widehat{P}^{(n)}_{n,n}$ of
the fundamental matrix. Then, according to
\[
     R^+_n = P_{n-1}\widehat{P}^{(n)}_{n,n},
\]
 we can have the equivalence. For details, readers may refer to \cite{[LR99]}, \cite{Zhao:2000} and \cite{[ZLB03]}.
\end{remark}}

{  The expression of the stationary distribution for the state-dependent walk in Theorem~\ref{thm3} enable us to obtain the following asymptotic behavior.}
Let $D=\{(P,Q,R):~(P+Q+R)\mathbf{1}=\mathbf{1},~\bar{\varrho}^{+}_1<\infty\}$.

\begin{theorem}\label{thm4}
{  For the random walk on a strip,}
\begin{description}
\item[(1)] If $(P,Q,R)\in D$, we have $\lambda_{A^{-}}<1$.

\item[(2)] Suppose {  that the} random walk starts from layer $0$ with {  the} censored measure $\mathbf{\breve{\mu}}_{0}$, and {  the transition probabilities satisfy} $(P_{n},Q_{n},R_{n})\rightarrow(P,Q,R)$ as $n\rightarrow\infty$ {  with} $(P,Q,R)\in D$. Then the random walk is positive recurrent and
the stationary distribution $\{\mathbf{\nu}_{n},~n\geq 0\}$ defined in (\ref{vn}) is light-tailed, with the decay rate {  $0<\lambda_{A^{-}} \leq 1$} along the layer direction, that is, for each fixed $1\leq j\leq d$,
\begin{equation}
\lim_{n\rightarrow\infty}\frac{\log\mathbf{\nu}_{n}(j)}{n}={  \log \lambda_{A^{-}},}
\end{equation}
where $\lambda_{A^{-}}$ is the maximum eigenvalues of $A^{-}$ (given in(\ref{a})).
\end{description}
\end{theorem}

{
\begin{example}[A retrial queue with a state-dependent retrial rate] \rm This model is a modification of the standard $M/M/c$ retrial queue (for example, see Falin and Templeton~\cite{Falin-Templeton:97}). In the modified model,
instead of the retrial rate $n \theta$, we assume the total retrial rate is $\theta_n$, where $n$ is the number of customers in the retrial orbit.
For this model, let $N(t)$ and $C(t)$ be the number of retrial customers in the orbit and the number of busy servers at time $t$, respectively.
Then, it is easy to see that $(N(t), C(t))$ is a continuous-time Markov chain. We show how to apply Theorem~\ref{thm4} to obtain the exponential decay rate.
For this purpose, assume that $\theta_n \to \theta<\infty$ as $n \to \infty$. Then, the generator of the limiting chain is given by
\begin{equation} \label{eqn:Q}
    Q = \left ( \begin{array}{ccccccc}
B_0 & A \\
C & B & A \\
& C & B & A \\
& & \ddots & \ddots & \ddots
\end{array} \right ),
\end{equation}
where
\[
    B = \left ( \begin{array}{ccccccc}
- (\lambda+\theta) & \lambda \\
\mu & -(\lambda+\mu+ \theta) & \lambda \\
 & \ddots & \ddots & \ddots \\
 & & (c-1) \mu & -[\lambda+(c-1)\mu+ \theta] & \lambda \\
 & & & c \mu & -(\lambda+c\mu)  \\
\end{array} \right ),
\]
\[
    A = \left ( \begin{array}{ccccccc}
0 &  \\
& 0 & \\
& & \ddots & \\
& & & 0 & \\
& & & & \lambda
\end{array} \right )
\qquad \mbox{and} \qquad
    C = \left ( \begin{array}{ccccccc}
0 &  \theta \\
& 0 &  \theta \\
& & \ddots & \ddots \\
& & & 0 &  \theta \\
& & & & 0
\end{array} \right ).
\]
Without loss of generality, we assume $\lambda+c \mu +\theta=1$. Upon uniformization, we can convert the generator to a transition matrix $\breve{P}=I-Q$ to have
$(P,Q,R)$. To determine the condition for positive recurrence and $\lambda_{A^-}$, we use Theorem~\ref{thm3} and (\ref{a}), respectively. First, (\ref{a}) is equivalent to the equation $R^+ = P + R^+R+{R^+}^2Q$. To solve this equation, we notice that
\[
    R^+=  \left [ \begin{array}{cccc}
0 & 0 & \cdots & 0 \\
\vdots & \vdots & \cdots & \vdots \\
0 & 0 & \cdots & 0 \\
r_1&r_2 &\cdots & r_c
\end{array} \right ],
\]
which greatly simplify the calculations. Also, we can find that the chain is positive recurrent if and only if $r_c<1$. For example, when $c=1$,
$\lambda_{A^-}=r_c =\lambda (\lambda+\theta)/\mu\theta$ and when $c=2$,
\[
    \lambda_{A^-}=r_c = \frac{\lambda}{\theta  \mu} \frac{(\lambda+\theta)^2+\theta \mu}{3\lambda +2\mu+2\theta}.
\]
 As $c$ gets larger, the formula becomes cumbersome and is less interesting.
\end{example}
}

We arrange the remainder of this paper as follows. As the main
tool of this paper, the intrinsic branching structure within
random walk on a strip is briefly reviewed in Section~\ref{s2};
and then the proofs for the Theorems are  followed in
{  Section}~\ref{s3}.

\section{A brief review for the intrinsic branching structure within random walk on a strip \label{s2}}
\setcounter{equation}{0}

The intrinsic branching structure within a random walk is a very powerful tool in the research {  on} the limit property about random walk. For the neighborhood nearest random walk on the line,  Dwass (\cite{[Dwa75]}, 1975) and Kesten {  et al.} (\cite{[KKS75]}, 1975) observed a Galton-Watson process {  with a} geometric offspring distribution hidden in it.
{  Kesten \textit{et al.}} (\cite{[KKS75]}) proved a stable law {  by using} the branching structure for the random walk {  in a} random environment.
{  For other random walks, e.g.,} the random walk with bounded jumps, the branching structure were revealed by Hong \& Wang (\cite{[HW09]},  2009) for the $(L,1)$-case and Hong \& Zhang (\cite{[HZ10]}, 2010) for the $(1,R)$-case.

{  The intrinsic} branching structure within random walk
on a strip has been revealed by Hong  \& Zhang (\cite{[HM12]},
2012), which {  enables us to provide explicit} criteria
for (positive) recurrence and {  to obtain an explicit
expression for the stationary distribution. As a consequence, it
allows us} to consider the tail asymptotic of the stationary
distribution. {  The key point} is the trajectory
decomposition for the random walk. If the walk starts from layer
$n>0$, the trajectory {  has the}  ``upper" and ``lower"
parts, {  which are introduced in the following
subsections.}

\subsection{The ``lower" branching structure}

Assume that $X_{0}\in L_{k}$, {  the} random walk starts from layer $k$ with initial distribution $\mathbf{\mu}_{k}$ {  or} $\mu_{k}(i)=P(\xi_{0}=k,~Y_{0}=i)$. For $0< n\leq k $ and $i\in \{1,2,\cdots,d\}$, define $U_{n}^{i}$ as the number of steps from layer $n$ to $(n-1,i)$ before {  the} hitting time $T_{k+1}$, {  and} $Z_{n}^{i}$ as the number of steps from layer $n$ to $(n,i)$ before $T_{k+1}$. Define
\begin{equation}
\mathbf{U}_{n}=(U_{n}^{1},U_{n}^{2},\cdots,U_{n}^{d}),\quad
\mathbf{Z}_{n}=(Z_{n}^{1},Z_{n}^{2},\cdots,Z_{n}^{d}),\no
\end{equation}
and $|\mathbf{U}_{n}|=\mathbf{U}_{n}\mathbf{1}$,
~$|\mathbf{Z}_{n}|=\mathbf{Z}_{n}\mathbf{1}$.

\begin{theorem}\label{thm11} (Hong $\&$ Zhang, 2012)~~
Suppose {  that} Condition~C is satisfied, {  and the} random walk starts from layer $k$ with initial distribution $\mathbf{\mu}_{k}$. Then
$\{|\mathbf{U}_{n}|,~1<n\leq k\}$ and
$\{|\mathbf{Z}_{n}|,~1<n\leq k\}$ are inhomogeneous branching processes with immigration. The offspring distribution ($1<n\leq k$) is given as:
\begin{eqnarray}
P\big(|\mathbf{U}_{n}|=m \big| \mathbf{U}_{n+1}=\mathbf{e}_{i}\big )&=&\mathbf{e}_{i} [(I-R_{n})^{-1}  Q_{n}\zeta_{n-1}^{+}]^{m}(I-R_{n})^{-1}P_{n}\mathbf{1},\no\\
P\big(|\mathbf{Z}_{n}|=K \big| \mathbf{U}_{n+1}=\mathbf{e}_{i}\big)&=&\mathbf{e}_{i}[(I-Q_{n}\zeta_{n-1}^{+})^{-1} R_{n}]^{K} (I-Q_{n}\zeta_{n-1}^{+})^{-1} P_{n} \mathbf{1},\no
\end{eqnarray}
with immigration
\begin{equation}
P\big(\mathbf{U}_{k+1}=\mathbf{e}_{i} \big )=\mu_{k}(i),\quad i\in \mathscr{D}.\no
\end{equation}
\qed
\end{theorem}

{  The key idea in the construction of the} branching mechanism is that the {  position of the walk corresponds} to the time of the branching process. $|\mathbf{U}_{n}|$ as the number of steps from layer $n$ to $n-1$ layer is indeed the $(k-n)$-th generation of the branching process. The condition $T_{k+1}<\infty$ is {   obviously} satisfied in our reflecting model.

\begin{proposition}\label{p22}~Denote $N_{n}^{i}$ as the number of steps visited {  at} $(n,i)$ before time $T_{k+1}$, and $\mathbf{N}_{n}=(N_{n}^{1},N_{n}^{2},\cdots,N_{n}^{d})$ {  with}
$|\mathbf{N}_{n}|=\mathbf{N}_{n}\mathbf{1}$. Suppose {  that} Condition~C is satisfied, {  and the} random walk starts from layer $k$ with initial distribution $\mathbf{\mu}_{k}$. Then for any $0<n\leq k$,
\begin{eqnarray}\label{nz}
E_{\mu_k}(\mathbf{N}_{n})&=&\mathbf{\mu}_{k}A_{k}^{+}A_{k-1}^{+}\cdots A_{n+2}^{+}A_{n+1}^{+}(I-Q_{n}\zeta_{n-1}^{+}-R_{n})^{-1},\no\\
~\no\\
E_{\mu_k}(|\mathbf{N}_{n}|)
&=&\mathbf{\mu}_{k}A_{k}^{+}A_{k-1}^{+}\cdots A_{n+2}^{+}A_{n+1}^{+}\mathbf{u}_{n}^{+},
\end{eqnarray}
and for $n=0$,
\begin{equation}\label{nz2}
E_{\mu_k}(\mathbf{N}_{0})=\mathbf{\mu}_{k}A_{k}^{+}A_{k-1}^{+}\cdots A_{1}^{+},\quad
E_{\mu_k}(|\mathbf{N}_{0}|)
=\mathbf{\mu}_{k}A_{k}^{+}A_{k-1}^{+}\cdots A_{1}^{+}\mathbf{1}.
\end{equation}
\end{proposition}

\proof {  We provide key steps here and readers may refer to} \cite{[HM12]} for details.
 Note that
\begin{eqnarray}\label{ma}
E(\mathbf{U}_{n}|\mathbf{U}_{n+1})
&=&\mathbf{U}_{n+1}\sum_{m=1}^{+\infty} [(I- R_{n})^{-1}Q_{n}\zeta_{n-1}^{+}]^{m-1}(I- R_{n})^{-1}Q_{n}
=\mathbf{U}_{n+1}A_{n}^{+},\\
E(\mathbf{Z}_{n}|\mathbf{U}_{n+1})
&=&\mathbf{U}_{n+1}\sum_{K=1}^{+\infty} [(I-Q_{n}\zeta_{n-1}^{+} )^{-1}
R_{n}]^{m-1} (I-Q_{n}\zeta_{n-1}^{+} )^{-1} R_{n}\no\\
&=&\mathbf{U}_{n+1}(I- Q_{n}\zeta_{n-1}^{+}-R_{n})^{-1} R_{n}.\no
\end{eqnarray}
With the help of {  the} branching structure in Theorem~\ref{thm11}, we have for any $0<n\leq k$,
\begin{eqnarray}
E_{\mu_k}(\mathbf{N}_{n})
&=& E_{\mu_k}(\mathbf{U}_{n}\zeta_{n-1}^{+}+\mathbf{Z}_{n}+\mathbf{U}_{n+1})\\
&=&E_{\mu_k}\big[E_{\mu_k}(\mathbf{U}_{n}\big|  \mathbf{U}_{n+1})\zeta_{n-1}^{+}+E_{\mu_k}(\mathbf{Z}_{n}\big|  \mathbf{U}_{n+1})+E_{\mu_k}(\mathbf{U}_{n+1}\big|\mathbf{U}_{n+1})\big]\no\\
&=&E_{\mu_k}(\mathbf{U}_{n+1})(I-Q_{n}\zeta_{n-1}^{+}-R_{n})^{-1},\no
\end{eqnarray}
\noindent and then $E_{\mu_k}(|\mathbf{N}_{n}|)=E_{\mu_k}(\mathbf{U}_{n+1})\mathbf{u}_{n}$.

For $n=0$, the expected number of steps visiting the reflecting layer $0$ before time $T_{k+1}$ is $E_{\mu_k}(\mathbf{N}_{0})=E_{\mu_k}(\mathbf{Z}_{0}+\mathbf{U}_{1})=E_{\mu_k}(\mathbf{U}_{1})$, and then $E_{\mu_k}(|\mathbf{N}_{0}|)=E_{\mu_k}(\mathbf{U}_{1})$. Together with
\begin{equation}
E_{\mu_k}(\mathbf{U}_{n+1})= E_{\mu_k}[E_{\mu_k}(\mathbf{U}_{n+1}| \mathbf{U}_{n+2})]=E_{\mu_k}(\mathbf{U}_{n+2})A_{n+1}^{+}=
\cdots=\mathbf{\mu}_{k}A_{k}^{+}A_{k-1}^{+}\cdots A_{n+2}^{+}A_{n+1}^{+},\no
\end{equation}
the proof is finished. We usually denote $\mathbf{u}_{0}^{+}=\mathbf{1}$ to make the expression uniform. \qed

\subsection{The ``upper" branching structure}
Assume that $X_{0}\in L_{k}$, {  the} random walk starts from layer $k$ with initial distribution $\mathbf{\mu}_{k}$ {  or} $\mu_{k}(i)=P(\xi_{0}=k,~Y_{0}=i)$.
Similarly, For $n\geq k+1$ and $i\in \{1,2,\cdots,d\}$, define
\begin{equation}
\mathbf{W}_{n}=(W_{n}^{1},W_{n}^{2},\cdots,W_{n}^{d}),\quad\mbox{and}\quad
\mathbf{Z}_{n}^{-}=(Z_{n}^{-,1},Z_{n}^{-,2},\cdots,Z_{-,n}^{d}),\no
\end{equation}
where $W_{n}^{i}$ is the number of steps from layer $n$ to $(n+1,i)$ before {  the} hitting time $T_{k-1}$; and $Z_{n}^{-,i}$ is the number of steps from layer $n$ to $(n,i)$ before {  the} hitting time $T_{k-1}$.

Define $|\mathbf{W}_{n}|=\mathbf{W}_{n}\mathbf{1}$ and
$|\mathbf{Z}_{n}^{-}|=\mathbf{Z}_{n}^{-}\mathbf{1}$.

\begin{theorem}\label{thm12} (Hong $\&$ Zhang, 2012)~~
Suppose {  that} Condition~C is satisfied, {  and the} random walk starts from layer $k$ with initial distribution $\mathbf{\mu}_{k}$ and $T_{k-1}<+\infty$. Then
$\{|\mathbf{W}_{n}|,~n\geq k+1\}$ and
$\{|\mathbf{Z}_{n}^{-}|,~n\geq k+1\}$ are inhomogeneous branching processes with immigration. The offspring distribution ($n\geq k+1$) is given by
\begin{eqnarray}
P\big(|\mathbf{W}_{n}|=m \big| \mathbf{W}_{n-1}=\mathbf{e}_{i}\big)
&=&\mathbf{e}_{i} [(I-R_{n})^{-1}  P_{n}\zeta_{n+1}^{-}]^{m}(I-R_{n})^{-1}Q_{n}\mathbf{1},\no\\
P\big(|\mathbf{Z}_{n}^{-}|=K \big| \mathbf{W}_{n-1}=\mathbf{e}_{i}\big)
&=&\mathbf{e}_{i}[(I-P_{n}\zeta_{n+1}^{-})^{-1} R_{n}]^{K} (I-P_{n}\zeta_{n+1}^{-})^{-1} Q_{n} \mathbf{1},\no
\end{eqnarray}
with immigration
\begin{equation}
P\big(\mathbf{U}_{k-1}=\mathbf{e}_{i}\big)=\mathbf{\mu}_{k}(i),\quad i\in \mathscr{D}.\no
\end{equation} \qed
\end{theorem}

{  In parallel,} for $n\geq k+1$, denote
\begin{equation}
N_{n}^{-,i}=\sharp\{ k \in [0,T_{-1}):~X_{k}=(n,i)\},\no
\end{equation}
where $N_{n}^{-,i}$ is the number of steps visited {  at} $(n,i)$ before {  the} hitting time $T_{k-1}$. Define $\mathbf{N}_{n}^{-}=(N_{n}^{-,1},N_{n}^{-,2},\cdots,N_{n}^{-,d})$ and $|\mathbf{N}_{n}^{-}|=\mathbf{N}_{n}^{-}\mathbf{1}$.

\begin{proposition}~
Suppose {  that} Condition~C is satisfied, {  and the} random walk starts from layer $k$ with initial distribution $\mathbf{\mu}_{k}$. Then for any $n\geq k+1$,
\begin{eqnarray}\label{nf}
E_{\mu_k}(\mathbf{N}_{n}^{-})
&=&E_{\mu_k}(\mathbf{W}_{n-1})(I-P_{n}\zeta_{n+1}^{-}-R_{n})^{-1}\no\\
&=&\mathbf{\mu}_{k}A_{k}^{-}A_{k+1}^{-}\cdots A_{n-2}^{-}A_{n-1}^{-}(I-P_{n}\zeta_{n+1}^{-}-R_{n})^{-1},\no\\
~\no\\
E_{\mu_k}(|\mathbf{N}_{n}^{-}|)
&=&E_{\mu_k}(\mathbf{W}_{n-1})(I-P_{n}\zeta_{n+1}^{-}-R_{n})^{-1}\mathbf{1}
=\mathbf{\mu}_{k}A_{k}^{-}A_{k+1}^{-}\cdots A_{n-2}^{-}A_{n-1}^{-}\mathbf{u}_{n}^{-}.\no\\
\end{eqnarray}
\end{proposition}~

\section{Proofs}\label{s3}
\setcounter{equation}{0}

\subsection{Criteria for recurrence---{\it Proof of Theorem \ref{thm1}}}

Let $T_{y}^{0}=0$, {  and let} $T_{y}^{k}=\inf\{n>T_{y}^{k-1}:~X_{n}\in L_{y}\}$ for $k\geq 1$, {  or} $T_{y}^{k}$ is the time of the $k$-th return to layer $y$. Note that $T_{y}^{1}>0$. {  Hence, a possible} visit at time $0$ does not count, {  and} $T_{y}^{1}$ equals to $T_{y}^{+}$ defined above.

We firstly extend a basic {  property} about Markov {  chains to the} random walk on a strip, which is stated {  in} the following lemma.
\begin{lemma}\label{lem1}
Layer $y$ is recurrent if and only if $E_{y}(|\mathbf{N}_{y}|)=+\infty$.
\end{lemma}

\proof  {  Denote $f_{x,y}=P_{x}(T_{y}^{+}<\infty)$. Then,}
\begin{equation}
P_{x}(T_{y}^{k}<\infty)=f_{x,y}f_{y,y}^{k-1}.\no
\end{equation}
{  This is clear, since in order to
visit layer $y$ for exactly the $k$-th time,} the walk has to go from layer $x$ to layer $y$ {  first, and then return to layer $y$ $k-1$ times. A detailed formal proof is similar to that for Theorem~3.1 in \cite{[Dur]} for the random walk on a line.}

Recall that $|\mathbf{N}_{y}|=\mathbf{N}_{y}\mathbf{1}=\sum_{m=1}^{+\infty}I_{\{X_{m}\in L_{y}\}}$ is the number of visits to layer $y$ at positive times. By the definition, {  layer $y$} is transient if and only if $f_{y,y}<1$. Suppose {  that} layer $y$ is transient, then
\begin{eqnarray}
E_{x}(|\mathbf{N}_{y}|)
&=&\sum_{k=1}^{+\infty}P_{x}(\mathbf{N}_{y}\mathbf{1}\geq k)
=\sum_{k=1}^{+\infty}P_{x}(T_{y}^{k}<\infty)\no\\
&=&\sum_{k=1}^{+\infty}f_{x,y}f_{y,y}^{k-1}=\frac{f_{x,y}}{1-f_{y,y}}<+\infty.\no
\end{eqnarray}
Thus, layer $y$ is recurrent if and only if $E_{y}(|\mathbf{N}_{y}|)=+\infty$.
\qed

{\it Proof of Theorem \ref{thm1}} ~ Because the random walk is irreducible, we {  only need} to calculate the $E_{\mu}(|\mathbf{N}_{0}|)$, where the walk {  starts} at layer 0 with distribution $\mu$, and $|\mathbf{N}_{0}|=\sum_{i=0}^{\infty}1_{(X_i\in L_0)}$ is the occupation time of the walk at layer $0$. We can decompose the trajectory of the walk as the summation of infinite ``pieces", each ``piece" is an immigration (``lower") branching structure as considered in Theorem \ref{thm11}. {  In fact, by recalling the definition of the hitting times $T_k=\inf\{i: X_i\in L_k\}$ for layer $k$ and denoting $X^{{  (\tau_k)}}=\{X_i, T_k< i \leq T_{k+1}\}$, we can write}
\begin{eqnarray}\label{31}
\{X_i, i> 0\}=\bigcup_{k=0}^{+\infty} \{X_i, T_k< i \leq T_{k+1}\}=\bigcup_{k=0}^{+\infty} X^{{  (\tau_k)}},
\end{eqnarray}
and as a consequence,
\begin{eqnarray}\label{32}
|\mathbf{N}_{0}|=\sum_{i=0}^{+\infty}1_{(X_i\in L_0)}=\sum_{k=0}^{+\infty}\sum_{i=T_k}^{T_{k+1}}1_{(X_i\in L_0)}=\sum_{k=0}^{+\infty}|\mathbf{N}_{0}|^{{  (\tau_k)}},
\end{eqnarray}
{  where the superscript is used to emphasize that the process starts at $T_k$.}
For $k=0,1$, {  \ldots,}  each {  trajectory ``piece"} $X^{{  (\tau_k)}}=\{X_i, T_k< i \leq T_{k+1}\}$ {  formulates} a branching structure with immigration $P\big(\mathbf{U}_{k+1}=\mathbf{e}_{i} \big )=\mu_{k}(i)$, where $\mathbf{\mu}_{k}=
\mathbf{\mu}\zeta_{0}^{-}\zeta_{1}^{+}\cdots\zeta_{k-1}^{+}$.  By (\ref{nz2}) of Proposition \ref{p22} we have,
$$E_{\mu_k}(|\mathbf{N}_{0}|^{{  (\tau_k)}})
=\mathbf{\mu}_{k}A_{k}^{+}A_{k-1}^{+}\cdots A_{1}^{+}\mathbf{1}.$$
Combining with (\ref{32}),
\begin{eqnarray}\label{33}
E_{\mu}|\mathbf{N}_{0}|=\sum_{k=0}^{+\infty}E_{\mu_k}|\mathbf{N}_{0}|^{{  (\tau_k)}}=\sum_{k=0}^{+\infty}\mathbf{\mu}_{k}A_{k}^{+}A_{k-1}^{+}\cdots A_{1}^{+}\mathbf{1}=\beta^{+}.\no
\end{eqnarray}
The proof is complete. \qed

\subsection{Criteria for positive recurrence---{\it Proof of Theorem \ref{thm2}}}

Define $\overline{T}_{n}^{n}$ as {  the return time} of layer $n$ when the random walk starting from layer $n$,
 $\overline{T}_{n}^{(n-1,j)}$  {  the hitting time} of layer $n$ when the random walk starting from $(n-1,j)$,
 {  and} $\overline{T}_{n}^{(n+1,j)}$  {  the hitting time} of layer $n$ when the random walk starting from $(n+1,j)$. Then {  by the path} decomposition,
\begin{equation}
\overline{T}_{n}^{n}=\sum_{j}I_{X_{1}=(n-1,j)}\overline{T}_{n}^{(n-1,j)}
+\sum_{j}I_{X_{1}=(n+1,j)}\overline{T}_{n}^{(n+1,j)}+\sum_{j}I_{X_{1}=(n,j)},\no
\end{equation}
and {  therefore,}
\begin{eqnarray}
E_{\mu_n}(T_{n}^{+})&=&E_{\mu_n}(\overline{T}_{n}^{n})\no\\
&=&E_{\mu_n}(\sum_{j}I_{X_{1}=(n-1,j)}(\overline{T}_{n}^{(n-1,j)}+1)
+\sum_{j}I_{X_{1}=(n+1,j)}(\overline{T}_{n}^{(n+1,j)}+1)+\sum_{j}I_{X_{1}=(n,j)})\no\\
&=& \sum_{j}P_{\mu_n}(X_{1}=(n-1,j))(E_{(n-1,j)}T_{n}^{+}+1)\no\\
&&+\sum_{j}P_{\mu_n}(X_{1}=(n+1,j))(E_{(n+1,j)}T_{n}^{+}+1)
+\sum_{j}P_{\mu_n}(X_{1}=(n,j)).\no
\end{eqnarray}
Note that
\begin{equation}
E_{(n-1,j)}(T_{n}^{+})
=\sum_{k=0}^{n-1}E_{(n-1,j)}(|\mathbf{N}_{k}|) \quad \mbox{and}\quad
E_{(n+1,j)}(T_{n}^{+})
=\sum_{k =n+1}^{+\infty}E_{(n+1,j)}(|\mathbf{N}_{n}^{-}|).\no
\end{equation}
It follows from (\ref{nz}) and (\ref{nf}) that
\begin{eqnarray}
E_{(n-1,j)}(T_{n}^{+})
&=&\mathbf{e}_{j}\sum_{k=0}^{n-1}A_{n-1}^{+}A_{n-2}^{+} \cdots A_{k+2}^{+}A_{k+1}^{+}\mathbf{u}_{k}^{+},\no\\
E_{(n+1,j)}(T_{n}^{+})
&=&\mathbf{e}_{j}\sum_{k =n+1}^{+\infty}A_{n+1}^{-}A_{n+2}^{-} \cdots A_{k-1}^{-}\mathbf{u}_{k}^{-}.\no
\end{eqnarray}
Hence  we have
\begin{eqnarray}
E_{\mu_n}(T_{n}^{+})
&=&\mathbf{\mu}_{n} Q_{n}(\sum_{k=0}^{n-1}
A_{n-1}^{+}A_{n-2}^{+} \cdots A_{k+1}^{+}\mathbf{u}_{k}^{+})
+\sum_{j}P_{\mu_n}\big(X_{1}=(n-1,j)\big)\no\\
&&+\mathbf{\mu}_{n}  P_{n}(\sum_{k =n+1}^{+\infty} A_{n+1}^{-}A_{n+2}^{-} \cdots A_{k-1}^{-}\mathbf{u}_{k}^{-})+\sum_{j}P_{\mu_n}\big(X_{1}=(n+1,j)\big)\no\\
&&+\sum_{j}P_{\mu_n}\big(X_{1}=(n,j)\big)\no\\
&=&\mathbf{\mu}_{n} Q_{n}(\sum_{k=0}^{n-1}
A_{n-1}^{+}A_{n-2}^{+} \cdots A_{k+1}^{+}\mathbf{u}_{k}^{+})
+\mathbf{\mu}_{n}  P_{n}(\sum_{k =n+1}^{+\infty} A_{n+1}^{-}A_{n+2}^{-} \cdots A_{k-1}^{-}\mathbf{u}_{k}^{-})
+\mathbf{\mu}_{n}\mathbf{1}.\no
\end{eqnarray}
{  Particularly, if the} random walk starts from layer $0$ with an initial distribution $\mathbf{\mu}$, we have
\begin{equation}\label{t0}
E_{\mu}(T_{0}^{+})
=\mathbf{\mu} P_{0}(\sum_{k \geq 1}A_{1}^{-}A_{2}^{-} \cdots A_{k-1}^{-}\mathbf{u}_{k}^{-})+\mathbf{\mu}\mathbf{1}.
\end{equation}
Thus the reflecting random walk on a strip is positive recurrent (independent of the initial distribution $\mu$) if and only if
$\varrho^{+}_1<\infty$ . \qed

\subsection{Stationary distribution---{\it Proof of Theorem \ref{thm3}}}

Suppose {  that the} random walk starts from layer $0$
with a censored measure $\mathbf{\breve{\mu}}_{0}$, which
satisfies
\begin{equation}
\mathbf{\breve{\mu}}_{0}\breve{P}^{(S_{0})}=\mathbf{\breve{\mu}}_{0}.\no
\end{equation}

The following lemma {  modifies} Thm~5.4.3 in \cite{[Dur]}
about {  a} stationary measure of {  a} general
Markov chain on $\mathbb{Z}^{d}$  to our model,
and defines a stationary measure for {  the} random walk
on a strip.
\begin{lemma}\label{lem2}
Suppose {  that the} random walk starts from layer $0$
with a censored measure $\mathbf{\breve{\mu}}_{0}$, and layer $0$
is a recurrent layer. Then $\{\mathbf{\overline{\nu}}_{n},~n\in
\mathbb{N}\}$ defines a stationary measure, where
\begin{equation}\label{nu}
\overline{\nu}_{n}(i)
=E_{\mathbf{\breve{\mu}}_{0}}\left(\sum_{m=0}^{T_{0}^{+}-1}I_{\{X_{m}=(n,i)\}}\right)
=\sum_{m=0}^{+\infty}P_{\mathbf{\breve{\mu}}_{0}}(X_{m}=(n,i),~m<T_{0}^{+}).
\end{equation}
\end{lemma}

\proof The {  key idea of the proof is to use the}
``cycle trick". $\overline{\nu}_{n}(i)$ is the expected number of
visits to $(n,i)$ {  at times $0,1,\ldots,T_{0}^{+}-1$.}
And $\sum_{y,j}\overline{\nu}_{y}(j)p\big[(y,j),(n,i)\big]$ is the
expected number of visits to $(n,i)$ {  at times $1, 2,
\ldots,T_{0}^{+}$,} which equals to $\overline{\nu}_{n}(i)$ since
$X_{T_{0}^{+}}\sim \mathbf{\breve{\mu}}_{0}$ if $X_{0}\sim
\mathbf{\breve{\mu}}_{0}$ {  based on the property of}
$\mathbf{\breve{\mu}}_{0}\breve{P}^{(S_{0})}=\mathbf{\breve{\mu}}_{0}$.

The goal is to {  prove that}
$\mathbf{\overline{\nu}}_{n}$ defined in (\ref{nu}) {  is
a stationary measure,} that is,
\begin{equation}\label{eq}
\sum_{y,j}\overline{\nu}_{y}(j)p\big[(y,j),(n,i)\big]
=\overline{\nu}_{n}(i).
\end{equation}
By Fibini's theorem, we get
\begin{equation}
\sum_{y,j}\overline{\nu}_{y}(j)p\big[(y,j),(n,i)\big]
=\sum_{m=0}^{+\infty}\sum_{y,j}P_{\mathbf{\breve{\mu}}_{0}}(X_{m}=(y,j),~m<T_{0}^{+})
p\big[(y,j),(n,i)\big].\no
\end{equation}

\noindent Case 1: $n\neq 0$. {  In this case, we have}
\begin{eqnarray}
& & \sum_{y,j}P_{\mathbf{\breve{\mu}}_{0}}(X_{m}=(y,j),~m<T_{0}^{+})p\big[(y,j),(n,i)\big]\no\\
&=&\sum_{y,j}P_{\mathbf{\breve{\mu}}_{0}}(X_{m}=(y,j),~m<T_{0}^{+},~X_{m+1}=(n,i))\no\\
&=&P_{\mathbf{\breve{\mu}}_{0}}(T_{0}^{+}>m+1,~X_{m+1}=(n,i)),\no
\end{eqnarray}
and then
\begin{eqnarray}
\sum_{y,j}\overline{\nu}_{y}(j)p\big[(y,j),(n,i)\big]
&=&\sum_{m=0}^{+\infty}\sum_{y,j}P_{\mathbf{\breve{\mu}}_{0}}(X_{m}=(y,j),~m<T_{0}^{+})
p\big[(y,j),(n,i)\big]\no\\
&=&\sum_{m=0}^{+\infty}P_{\mathbf{\breve{\mu}}_{0}}\big(T_{0}^{+}>m+1,~X_{m+1}=(n,i)\big).\no\\
&=& \sum_{m=0}^{+\infty}P_{\mathbf{\breve{\mu}}_{0}}(X_{m}=(n,i),~m<T_{0}^{+})=\overline{\nu}_{n}(i),
\end{eqnarray}
because $P_{\mathbf{\breve{\mu}}_{0}}(T_{0}^{+}>0,~X_{0}=(n,i))=0$.
\medskip

\noindent Case 2: $n=0$. At first, note that {  the
process starts} from layer 0 with the initial distribution
$\mathbf{\breve{\mu}}_{0}$, i.e., the right hand {  side}
of (\ref{eq}) is $\mathbf{\breve{\mu}}_{0}$. {  For the
left hand side} of (\ref{eq}), we calculate
\begin{eqnarray}
& &\sum_{y,j}P_{\mathbf{\breve{\mu}}_{0}}\big(X_{m}=(y,j),~m<T_{0}^{+}\big)p\big[(y,j),(n,i)\big]\no\\
&=&\sum_{y,j}P_{\mathbf{\breve{\mu}}_{0}}\big(X_{m}=(y,j),~m<T_{0}^{+},~X_{m+1}=(n,i)\big)\no\\
&=&P_{\mathbf{\breve{\mu}}_{0}}\big(T_{0}^{+}=m+1,~X_{m+1}=(0,i)\big),\no
\end{eqnarray}
and then
\begin{eqnarray}
\sum_{y,j}\overline{\nu}_{y}(j)p\big[(y,j),(n,i)\big]
&=&\sum_{m=0}^{+\infty}\sum_{y,j}P_{\mathbf{\breve{\mu}}_{0}}(X_{m}=(y,j),~m<T_{0}^{+})
p\big[(y,j),(n,i)\big]\no\\
&=&\sum_{m=0}^{+\infty}P_{\mathbf{\breve{\mu}}_{0}}(T_{0}^{+}=m+1,~X_{m+1}=(0,i)).\no
\end{eqnarray}

Note that $T_{0}^{+}\geq1$. {  Therefore,
$P_{\mathbf{\breve{\mu}}_{0}}(T_{0}^{+}=0,~X_{0}=(0,i))=0$, and}
we have
\begin{equation}
\sum_{m=0}^{+\infty}P_{\mathbf{\breve{\mu}}_{0}}\big(T_{0}^{+}=m+1,~X_{m+1}=(0,i)\big)
=\sum_{m=0}^{+\infty}P_{\mathbf{\breve{\mu}}_{0}}\big(T_{0}^{+}=m,~X_{m}=(0,i)\big)=\mathbf{\breve{\mu}}_{0},\no
\end{equation}
because $\mathbf{\breve{\mu}}_{0}$ is the censored measure.
{  The proof is complete now.} \qed

{  {\it Proof of Theorem \ref{thm3}}}~
{  First,} we can calculate the stationary measure in
(\ref{nu}) by {  using} the branching structure. The
stationary measure {  is given by}
$\overline{\nu}_{n}(i)=E_{\mathbf{\breve{\mu}}_{0}}\left(\sum_{m=0}^{T_{0}^{+}-1}I_{\{X_{m}=(n,i)\}}\right)$,
which is the expected number of visits to $(n,i)$ before time
$T_{0}^{+}$ (but not contains the time $T_{0}^{+}$). So,
$\overline{\mathbf{\nu}}_{n}(i)$~ $(n>0)$ equals to
$E_{1}\mathbf{N}_{n}^{-}$ obtained by {  the} branching
structure in (\ref{nf}). The stationary measure
$\{\mathbf{\bar{\mathbf{\nu}}}_{n},~n\in \mathbb{Z}\}$ can be
expressed as
\begin{equation}
\mathbf{\bar{\mathbf{\nu}}}_{n}=
\begin{cases}
\mathbf{\breve{\mu}}_{0}P_{0}A_{1}^{-}A_{2}^{-}\cdots A_{n-1}^{-}(I-P_{n}\zeta_{n+1}^{-}-R_{n})^{-1} & n> 0,\\
\mathbf{\breve{\mu}}_{0} & n= 0.
\end{cases}\no
\end{equation}
Note that
\begin{equation}
\sum_{n,i}\overline{\mathbf{\nu}}_{n}(i)
=\sum_{m=0}^{+\infty}P_{\mathbf{\breve{\mu}}_{0}}(T_{0}>m)=E_{\mathbf{\breve{\mu}}_{0}}T_{0}^{+}.\no
\end{equation}
The condition $\varrho^{+}_1<\infty$ {  ensures} that
\begin{equation}
\sum_{n,i}\overline{\mathbf{\nu}}_{n}(i)=E_{0}(T_{0}^{+})=
\mathbf{\breve{\mu}}_{0}P_{0}
\left(\sum_{k\geq1}A_{1}^{-}A_{2}^{-}\cdots
A_{k-1}^{-}\mathbf{u}_{k}^{-}\right)+\mathbf{\breve{\mu}}_{0}\mathbf{1}=\varrho^{+}<\varrho^{+}_1<\infty.\no
\end{equation}
{  As} a consequence the  stationary distribution
{  equals}
\begin{equation}
\mathbf{\nu}_{n}(i)=\frac{\overline{\nu}_{n}(i)}{E_{0}T_{0}^{+}}
=\frac{\mathbf{\breve{\mu}}_{0}P_{0}A_{1}^{-}A_{2}^{-}\cdots A_{n-1}^{-}\widetilde{u}_{n}^{-}}{\mathbf{\breve{\mu}}_{0} P_{0}(\sum_{k \geq 1}A_{1}^{-}A_{2}^{-} \cdots A_{k-1}^{-}\mathbf{u}_{k}^{-})+\mathbf{\breve{\mu}}_{0}\mathbf{1}}\quad n> 0,\no
\end{equation}
where $\widetilde{u}_{n}^{-}=(I-P_{n}\zeta_{n+1}^{-}-R_{n})^{-1}$, and
\begin{equation}
\mathbf{\nu}_{0}(i)=\frac{\mathbf{\breve{\mu}}_{0}}{\mathbf{\breve{\mu}}_{0} P_{0}(\sum_{k \geq 1}A_{1}^{-}A_{2}^{-} \cdots A_{k-1}^{-}\mathbf{u}_{k}^{-})+\mathbf{\breve{\mu}}_{0}\mathbf{1}}.\no
\end{equation}
\qed

\subsection{Light-tailed behavior---{\it Proof of Theorem \ref{thm4}}}

{  It is well-known that the stationary distribution for
the state-independent random walk (or a QBD process) on a
half-strip is matrix-geometric. Therefore, the tail has a
geometric (or exponential) decay. For the state-dependent random
walk on a half-strip, the stationary tail does not always have an
exponential decay. In this paper, we provide a criterion for this
case, which is proved here.}

\subsubsection{{  Preliminaries}}
Let $B=(b_{i,j})>0$ (which is {  called a} positive
matrix) if all $b_{i,j}>0$; and $B\geq0$ if all $b_{i,j}\geq 0$.
The spectrum of $n\times n$ matrix $B$ is denoted as
$\sigma(B)=\{\lambda_{1},\lambda_{2},\cdots \lambda_{n}\}$, where
$\sigma(B)$ is the set of all eigenvalues
$\lambda_{i}\in\mathbb{C}$. Define the spectral radius of $B$ as
$\rho(B)=\max\{|\lambda_{i}|:~\lambda_{i}\in\sigma(B),~1\leq i\leq
n\}$.

\begin{proposition} \label{prop1}({\it {  Perron's} Theorem} in \cite{[HJ90]})~If {  $B>0$ is an $n\times n$ matrix}, then
\begin{description}
\item[(1)] $\rho(B)>0$;

\item[(2)] $\rho(B)$ is an eigenvalue of $B$, and it is the unique
eigenvalue of maximum modulus;

\item[(3)] $\rho(B)$ is {   algebraically (and hence
geometrically) simple;}

\item[(4)]
\begin{equation}
\lim_{m\rightarrow\infty}\left(\frac{B}{\rho(B)}\right
)^{m}={  L>0.}\no
\end{equation}
\end{description}
\end{proposition}

Define $\|\cdot\|$ as the maximum column sum matrix norm, i.e. $\|B\|=\max_{1\leq i\leq n}\sum_{j=1}^{n}|b_{i,j}|$ for $B=(b_{i,j})$.

\begin{proposition}\label{prop2} (Krause (\cite{[Kra94]}, 94), Ostrowski (\cite{[Ost73]}, 73)) ~
Denote $\sigma(A)=\{\lambda_{1},\lambda_{2},\cdots \lambda_{n}\}$
where $\lambda_{i}$ are eigenvalues of $A$, {  and}
$\sigma(B)=\{\mu_{1},\mu_{2},\cdots \mu_{n}\}$,  where $\mu_{i}$
are eigenvalues of $B$. Define $d(\sigma(A),\sigma(B))$ as
{  the optimal matching distance} between the
{  spectrums} $\sigma(A)$ and  $\sigma(B)$, that is,
\begin{equation}
d(\sigma(A),\sigma(B))=\min_{\theta \in S_{n}}\max_{1\leq i\leq n}|\lambda_{i}-\mu_{\theta_{i}}|,\no
\end{equation}
where $S_{n}$ is denoted as the group of all permutations on sets
$\{1,2,\cdots, n\}$. Then, for any two matrices $A, B\in
\mathbb{R}^{n\times n}$, we have
\begin{equation}
d(\sigma(A),\sigma(B))\leq 4(2K)^{1-\frac{1}{n}}\|A-B\|^{\frac{1}{n}},\no
\end{equation}
where $K=\max\{\|A\|,\|B\|\}$.
\end{proposition}

Simply speaking, Proposition \ref{prop2} tells us that there
exists a permutation $\theta \in S_{n}$, such that the maximum
{  distance between the corresponding eigenvalues} is
small enough.

\subsubsection{Spectral radius}

Consider {  the state-independent} random walk
$\{\overline{X}_{n},~n\geq 0\}$ with transition probability
{  block} $(P,Q,R)$, starting from layer $0$ with an
initial distribution $\mathbf{\bar{\mu}}$. Let $\zeta^{-}$ be the
unique sequence of stochastic matrices {  satisfying}
\begin{equation}
\zeta^{-}=(I-P\zeta^{-}-R)^{-1}Q,\no
\end{equation}
and
\begin{equation}
A^{-}=(I-P\zeta^{-}-R)^{-1}P,\quad\quad
\mathbf{u}^{-}=(I-P\zeta^{-}-R)^{-1}\mathbf{1}.\no
\end{equation}

Denote the spectral radius of $A^{-}$ as $\rho(A^{-})$, and the
maximum eigenvalues of $A^{-}$ as $\lambda_{A^{-}}$. Assume
{  that the} random walk is positive recurrent
($\bar{\varrho}^{+}_1<\infty$).

{\it Proof of (2) in Corollary \ref{cor5}} ~ By {\it
{  Perron's} Theorem} in Proposition \ref{prop1}, we have
$\rho(A^{-})=\lambda_{A^{-}}$. The condition says
 $\bar{\varrho}^{+}_1=  \mathbf{1}'P~(\sum_{k\geq1}~(A^{-})^{k-1}~\mathbf{u}^{-})
+d<\infty$, i.e.,
\begin{equation}\label{sg}
\mathbf{1}'P\cdot\left(\sum_{k\geq1}(\lambda_{A^{-}})^{k-1}
(\frac{A^{-}}{\lambda_{A^{-}}})^{k-1}\right)\cdot\mathbf{u}^{-}<\infty.
\end{equation}
{  On the other hand, from (4)} of
Proposition~\ref{prop1}, we know
\begin{equation}\label{lad}
\lim_{k\rightarrow\infty}
(\frac{A^{-}}{\lambda_{A^{-}}})^{k-1} =L>0
\end{equation}
{  which, together with (\ref{sg}), leads to}
$\lambda_{A^{-}}<1.$ \qed

Let
\begin{equation}
E=\left(
        \begin{array}{cccc}
               1 &1 & \cdots & 1\\
               1 &1 & \cdots & 1\\
               ~ &~ & \ddots& ~\\
               1 &1 & \cdots & 1\\
             \end{array}
           \right)_{d\times d}.\no
\end{equation}

Denote the maximum eigenvalues of {  $(A^{-}-\varepsilon
E)$} as $\lambda_{\varepsilon}^{-}$, and the maximum eigenvalues
of {  $(A^{-}+\varepsilon E)$} as
$\lambda_{\varepsilon}^{+}$. We {  then} have
\begin{equation}
\rho(A^{-}-\varepsilon E)=\lambda_{\varepsilon}^{-}\quad\mbox{and}\quad \rho(A^{-}+\varepsilon E)=\lambda_{\varepsilon}^{+}.\no
\end{equation}

\begin{lemma}\label{lem6}
Suppose {  that $\lambda_{A^{-}}<1$, and for
$A^{-}=(a^{-}_{i,j})$,
\[
    \varepsilon<\min\{\displaystyle\min_{i,j}a^{-}_{i,j},~\frac{1}{C^{d}}(1-\lambda_{A^{-}})^{d}\}.
\]
Let $C=4(2\|A^{-}+E\|)^{1-\frac{1}{d}} d^{\frac{1}{d}}$.}
 Then,
\begin{equation}\label{lam}
\lambda_{A^{-}}-C \varepsilon^{\frac{1}{d}}<\lambda_{\varepsilon}^{-}<1,
\quad\mbox{and}\quad
\lambda_{\varepsilon}^{+}<\lambda_{A^{-}}+C \varepsilon^{\frac{1}{d}}<1.
\end{equation}
\end{lemma}

\proof For such $\varepsilon>0$, both $A^{-}-\varepsilon E$ and
$A^{-}+\varepsilon E$ are {   positive matrices, and} the
definitions of $\lambda_{\varepsilon}^{-}$ and
$\lambda_{\varepsilon}^{+}$ are meaningful. By Proposition
\ref{prop1}, $\lambda_{\varepsilon}^{-}$ and
$\lambda_{\varepsilon}^{+}$ are real-valued. By Proposition
\ref{prop2}, {  it is} not hard to get that
\begin{equation}
d(\sigma(A^{-}),\sigma(A^{-}-\varepsilon E))\leq 4(2K)^{1-\frac{1}{d}} d^{\frac{1}{d}} \varepsilon^{\frac{1}{d}}\leq C \varepsilon^{\frac{1}{d}},\no
\end{equation}
and
\begin{equation}
d(\sigma(A^{-}),\sigma(A^{-}+\varepsilon E))\leq 4(2K)^{1-\frac{1}{d}} d^{\frac{1}{d}} \varepsilon^{\frac{1}{d}}\leq C \varepsilon^{\frac{1}{d}},\no
\end{equation}
where $C=4(2K)^{1-\frac{1}{d}} d^{\frac{1}{d}}$, $K=\|A^{-}+E\|$.

Note that $0 < A^{-}-\varepsilon E \leq A^{-} \leq A^{-}+\varepsilon E$, then
\begin{equation}
1>\lambda_{A^{-}}=\rho(A^{-})\geq \rho(A^{-}-\varepsilon E)=\lambda_{\varepsilon}^{-}\quad\mbox{and}\quad
\lambda_{A^{-}}=\rho(A^{-})\leq \rho(A^{-}+\varepsilon E)=\lambda_{\varepsilon}^{+},\no
\end{equation}
{|color{red}and} $\lambda_{A^{-}}+C \varepsilon^{\frac{1}{d}}<1$
for such $\varepsilon>0$.

{  It is obvious} that $\lambda_{A^{-}}-C
\varepsilon^{\frac{1}{d}}<\lambda_{\varepsilon}^{-}<1$.
{  Otherwise} if
$\lambda_{\varepsilon}^{-}<\lambda_{A^{-}}-C\varepsilon^{\frac{1}{d}}$,
then for $\lambda_{A^{-}}$, there exists no permutation such that
$d(\sigma(A^{-}),\sigma(A^{-}-\varepsilon E))\leq C
\varepsilon^{\frac{1}{d}}$, due to {  the fact that}
$\lambda_{\varepsilon}^{-}$ is the largest eigenvalue of
$A^{-}-\varepsilon E$.

Similarly, $\lambda_{\varepsilon}^{+}<\lambda_{A^{-}}+C
\varepsilon^{\frac{1}{d}}<1$ holds. {  Otherwise} if
$\lambda_{\varepsilon}^{+}>\lambda_{A^{-}}+C\varepsilon^{\frac{1}{d}}$,
then for $\lambda_{A^{-}}^{\varepsilon+}$, there exists no
permutation such that $d(\sigma(A^{-}),\sigma(A^{-}-\varepsilon
E))\leq C \varepsilon^{\frac{1}{d}}$ holds, due to {  the
fact that} $\lambda_{A^{-}}$ is the largest eigenvalue of $A^{-}$.
\qed

\subsubsection{Light-tailed behavior--- proof of Theorem \ref{thm4}}

Recall $D=\{(P,Q,R):~(P+Q+R)\mathbf{1}=\mathbf{1},~\bar{\varrho}^{+}_1<\infty\}$.

{\it Proof of Theorem \ref{thm4}}~ {  Part (1) of the
theorem has been} proved in  (2) of  Corollary~\ref{cor5}. Now, we
focus on {  part (2).} The random walk $\{X_{n},~n\in
\mathbb{Z}\}$ starts from layer $0$ with an censored measure
$\mathbf{\breve{\mu}}_{0}$, and the transition
{  probabilities: $(P_{n},Q_{n},R_{n})\rightarrow(P,Q,R)
\in D$ as $n\rightarrow\infty$. It} is easy to find that the
random walk $\{X_{n},~n\in \mathbb{Z}\}$ is positive recurrent. To
this end, recall that
$A_{n}^{-}=(I-P_{n}\zeta_{n+1}^{-}-R_{n})^{-1}P_{n}$ and
$A^{-}=(I-P\zeta^{-}-R)^{-1}P$. {  Then,}
$A_{n}^{-}\rightarrow A^{-}$ as $n\rightarrow +\infty$ because
that $(P_{n},Q_{n},R_{n})\rightarrow(P,Q,R)$ as
$n\rightarrow\infty$; and ${\varrho}^{+}_1<\infty$ follows from $
~\bar{\varrho}^{+}_1<\infty$ as $(P,Q,R)\in D$.

Also we have $\lambda_{A^{-}}<1$ as $(P,Q,R)\in D$. For each
$\varepsilon>0$ defined in Lemma \ref{lem6}, there exists $N$,
{  such that} when $n>N$,
\begin{equation}
0<A^{-}-\varepsilon E \leq A_{n}^{-} \leq A^{-}+\varepsilon E,\no
\end{equation}
and then
\begin{equation}
(A^{-}-\varepsilon E)^{k}
\leq A_{N+1}^{-}A_{N+2}^{-}\cdots A_{N+k}^{-}
\leq (A^{-}+\varepsilon E)^{k}.\no
\end{equation}
Let
\begin{equation}\label{phi}
\Phi_{n}(i)=\mathbf{\mu}P_{0}A_{1}^{-}A_{2}^{-}\cdots A_{n-1}^{-}\widetilde{u}_{n}^{-}(i).
\end{equation}
Now we consider the {  first inequality.} Notice
\begin{equation}
A_{N+1}^{-}A_{N+2}^{-}\cdots A_{N+k}^{-}
\geq (\lambda_{\varepsilon}^{-})^{k}(\frac{A^{-}-\varepsilon E}{\lambda_{\varepsilon}^{-}})^{k},\no
\end{equation}
for the given $N$, {  therefore} we have
\begin{eqnarray}
\Phi_{N+k}(i)
&=&\mathbf{\mu}P_{0}A_{1}^{-}A_{2}^{-}\cdots A_{N+k-1}^{-}\mathbf{u}_{N+k}^{-}
=\mathbf{\mu} D_{1}(N)
\cdot \big(A_{N+1}^{-}A_{2}^{-}\cdots A_{k+N-1}^{-}\big) \cdot \widetilde{u}_{N+k}^{-}(i)\no\\
&\geq & \mathbf{\mu} D_{1}(N)
(\lambda^{-}_{\varepsilon})^{k}(\frac{A^{-}-\varepsilon E}{\lambda^{-}_{\varepsilon}})^{k}~\widetilde{u}_{N+k}^{-}(i),\no
\end{eqnarray}
where $D_{1}(N)=P_{0}A_{1}^{-}A_{2}^{-}\cdots A_{N}^{-}$ and
$\widetilde{u}_{N+k}^{-}=(I-P_{N+k}\zeta_{N+k+1}^{-}-R_{N+k})^{-1}$.
Hence,
\begin{equation}\label{nk}
\frac{\log\Phi_{N+k}(i)}{N+k}
\geq \frac{k \log \lambda^{-}_{\varepsilon}}{N+k}
+\frac{\log \mathbf{\mu} D_{1}(N)(\frac{A^{-}-\varepsilon E}{\lambda^{-}_{\varepsilon}})^{k}~\widetilde{u}_{N+k}^{-}(i) }{N+k}.
\end{equation}

By Proposition \ref{prop1}, there exists a positive matrix $W_{\varepsilon}^{-}$, such that
\begin{equation}
\lim_{k\rightarrow\infty}(\frac{A^{-}-\varepsilon E}{\lambda_{\varepsilon}^{-}})^{k}=W_{\varepsilon}^{-}.\no
\end{equation}

Together with $\displaystyle\lim_{k\rightarrow\infty}\widetilde{u}_{N+k}^{-}(i)
=\widetilde{u}^{-}(i)=(I-P\zeta^{-}-R)^{-1}(i)$, as a consequence,

\noindent $\log \mathbf{\mu} D_{1}(N)(\frac{A^{-}-\varepsilon E}{\lambda^{-}_{\varepsilon}})^{k}~\widetilde{u}_{N+k}^{-}(i)$ is bounded in $k$. Thus from (\ref{nk}),
\begin{equation}
\liminf_{k\rightarrow\infty}\frac{\log\Phi_{N+k}(i)}{N+k}
\geq \log \lambda^{-}_{\varepsilon}.\no
\end{equation}

Note that from Lemma \ref{lem6}, $\lambda_{A^{-}}-C \varepsilon^{\frac{1}{d}}<\lambda_{\varepsilon}^{-}<1$, so
\begin{equation}\label{inf}
\liminf_{k\rightarrow\infty}\frac{\log\Phi_{N+k}(i)}{N+k}
\geq \log (\lambda_{A^{-}}-C \varepsilon^{\frac{1}{d}}).
\end{equation}

Similarly, for the {  second inequality,} notice
\begin{equation}
A_{N+1}^{-}A_{N+2}^{-}\cdots A_{N+k}^{-}
\leq (\lambda_{\varepsilon}^{+})^{k}(\frac{A^{-}-\varepsilon E}{\lambda_{\varepsilon}^{+}})^{k},\no
\end{equation}
{  therefore} we have
\begin{equation}
\limsup_{k\rightarrow\infty}\frac{\log\Phi_{N+k}(i)}{N+k}
\leq \log \lambda^{+}_{\varepsilon}.\no
\end{equation}
Note that from Lemma \ref{lem6}, $\lambda_{\varepsilon}^{+}<\lambda_{A^{-}}+C \varepsilon^{\frac{1}{d}}<1$, so
\begin{equation}\label{sup}
\limsup_{k\rightarrow\infty}\frac{\log\Phi_{N+k}(i)}{N+k}
\leq \log (\lambda_{A^{-}}+C \varepsilon^{\frac{1}{d}}).
\end{equation}

Combine (\ref{inf}) and (\ref{sup}) {  to have}
\begin{equation}
\log (\lambda_{A^{-}}-C \varepsilon^{\frac{1}{d}})
\leq\liminf_{k\rightarrow\infty}\frac{\log\Phi_{k}(i)}{k}
\leq\limsup_{k\rightarrow\infty}\frac{\log\Phi_{k}(i)}{k}
\leq \log(\lambda_{A^{-}}+C \varepsilon^{\frac{1}{d}}).\no
\end{equation}

Let $\varepsilon\rightarrow 0$, we get
\begin{equation}
\lim_{k\rightarrow\infty}\frac{\log\Phi_{k}(i)}{k}
=\log \lambda_{A^{-}}.\no
\end{equation}

If $(P_{n},Q_{n},R_{n})\rightarrow (P,Q,R)$,~$(P,Q,R)\in D$, the
stationary distribution $\{\mathbf{\nu}_{n},~n\geq 0\}$
{  is given by}
\begin{equation}
\mathbf{\nu}_{n}
=\frac{\mathbf{\breve{\mu}}_{0}P_{0}A_{1}^{-}A_{2}^{-}\cdots
A_{n-1}^{-}\widetilde{u}_{n}^{-}}{\mathbf{\breve{\mu}}_{0}P_{0}(\sum_{k
\geq 1}A_{1}^{-}A_{2}^{-} \cdots
A_{k-1}^{-}\mathbf{u}_{k}^{-})+\mathbf{\breve{\mu}}_{0}\mathbf{1}},\quad
n>0,\no
\end{equation}
where the denominator $\mathbf{\breve{\mu}}_{0} P_{0}(\sum_{k \geq 1}A_{1}^{-}A_{2}^{-} \cdots A_{k-1}^{-}\mathbf{u}_{k}^{-})+\mathbf{\breve{\mu}}_{0}\mathbf{1}<{\varrho}^{+}_1<+\infty$.

Thus for the stationary distribution $\{\mathbf{\nu}_{n},~n\geq 0\}$, we have
\begin{equation}
\lim_{n\rightarrow\infty}\frac{\log\mathbf{\nu}_{n}(i)}{n}
=\lim_{k\rightarrow\infty}\frac{\log\Phi_{k}(i)}{k}
=\log \lambda_{A^{-}},\no
\end{equation}
i.e., the stationary distribution is light-tailed, with the decay rate $0\leq\lambda_{A^{-}}<1$ along the layer direction. \qed


\end{document}